\documentclass[12pt, a4paper]{article}
\usepackage{amsmath,amssymb,amsbsy,amsfonts,amsthm,latexsym,
                          amsopn,amstext,amsxtra,euscript,amscd}

\begin{document}

\newtheorem{lem}{Lemma}
\newtheorem{lemma}[lem]{Lemma}
\newtheorem{prop}{Proposition}
\newtheorem{thm}{Theorem}
\newtheorem{theorem}[thm]{Theorem}

\def\Ep{\,\, \substack{+\\ E}\,\,}

\def\En{\,\, \substack{-\\ E}\,\,}

\def\Et{\,\, \substack{*\\ E}\,\,}

\def\Gn{\,\, \substack{-\\ G}\,\,}

\def\Gt{\,\, \substack{*\\ G}\,\,}

\title{\sc On the size of the set $A(A+1)$}

\author{M. Z. Garaev and Ch.-Y. Shen}

\author{
{\sc Moubariz Z. Garaev} \\
{Instituto de Matem{\'a}ticas}\\
{Universidad Nacional Aut\'onoma de M{\'e}xico} \\
{C.P. 58089, Morelia, Michoac{\'a}n, M{\'e}xico} \\
{\tt garaev@matmor.unam.mx}\\
\and
{\sc Chun-Yen Shen} \\
{Department of Mathematics}\\
{Indiana University} \\
{Rawles Hall, 831 East Third St.}\\
{Bloomington, IN 47405, USA}\\
{\tt shenc@indiana.edu}}

\maketitle

\begin{abstract}
Let $F_p$ be the field of a prime order $p.$ For a subset $A\subset
F_p$ we consider the product set $A(A+1).$ This set is an image of
$A\times A$ under the polynomial mapping $f(x,y)=xy+x:F_p\times
F_p\to F_p.$ In the present note we show that if $|A|<p^{1/2},$ then
$$
|A(A+1)|\ge |A|^{106/105+o(1)}.
$$
If $|A|>p^{2/3},$ then we prove that
$$
|A(A+1)|\gg \sqrt{p\, |A|}
$$
and show that this is the optimal in general settings bound up to
the implied constant. We also estimate the cardinality of $A(A+1)$
when $A$ is a subset of real numbers. We show that in this case one
has the Elekes type bound
$$
|A(A+1)|\gg |A|^{5/4}.
$$
\end{abstract}

\section {Introduction}

Let $F_p$ be the field of residue classes modulo a prime number $p$
and let $A$ be a non-empty subset of $F_p.$ It is known
from~\cite{BGK, BKT} that if $|A|<p^{1-\delta},$ where $\delta>0,$
then one has the sum-product estimate
$$
|A+A|+|AA|\gg |A|^{1+\varepsilon}; \qquad
\varepsilon=\varepsilon(\delta)>0.
$$
This estimate and its proof consequently have been quantified and
simplified in~\cite{BG},~\cite{Gar1}--\cite{HIS},~\cite{KS}--\cite{Sh},~\cite{TV}. From the
sum-product estimate and Ruzsa's triangle inequalities (see,~\cite{R1} and~\cite{R2}) it follows
that the polynomial $f(x,y,z)=xy+z:F_p^3\to F_p$ possesses an
expanding property, in a sense that for any subsets $A,B,C$ with
$|A|\sim |B|\sim |C|\sim p^{\alpha},$ where $0<\alpha<1$ is fixed,
the set $f(A,B,C)$ has cardinality greater than $p^{\beta}$ for
some $\beta=\beta(\alpha)>\alpha.$ The problem raised by Widgerson
asks to explicitly write a polynomial with two variables which
would satisfy the expanding condition. This problem was solved by
Bourgain~\cite{B1}, showing that one can take $f(x,y)=x^2+xy.$

Now consider the polynomial $f(x,y)=xy+x.$ This polynomial, of
course, does not possess the expanding property in the way defined
above. Nevertheless, from Bourgain's work~\cite{B1} it is known that
if $|A|\sim p^{\alpha},$ where $0<\alpha<1,$ then
$$
|f(A,A)|=|A(A+1)|\ge p^{\beta};\qquad \beta=\beta(\alpha)>\alpha.
$$

In the present note we deal with explicit lower bounds for the size
of the set $A(A+1).$ Our first result addresses the most nontrivial
case $|A|<p^{1/2}.$

\begin{theorem}
\label{thm:106/105} Let $A\subset F_p$ with $|A|<p^{1/2}.$ Then
$$
|A(A+1)| \ge |A|^{106/105+o(1)}.
$$
\end{theorem}

Theorem~\ref{thm:106/105} will be derived from the
Balog-Szemer\'edi-Gowers type estimate and a version of the
sum-product estimate given in~\cite{BG}. We remark that the
statement of Theorem~\ref{thm:106/105} remains true in a slightly
wider range than $|A|<p^{1/2}.$ On the other hand, if $|A|>p^{2/3},$ then we
have the optimal in general settings bound.

\begin{theorem}
\label{thm:optimal} For any subsets $A, B, C\subset F_p^*$ the
following bound holds:
$$
|AB|\cdot|(A+1)C|\gg \min\Bigl\{p\,|A|,\,
\frac{|A|^2\cdot|B|\cdot|C|}{p}\Bigr\}.
$$
\end{theorem}

Theorem~\ref{thm:optimal} can be compared with the following
estimate from~\cite{Gar2}:
$$
|A+B|\cdot |AC|\gg \min\Bigl\{p\, |A|,\,
\frac{|A|^2\cdot|B|\cdot|C|}{p}\Bigr\}.
$$

Taking $B=A+1,\, C=A,$ Theorem~\ref{thm:optimal} implies
$$
|A(A+1)|\gg \min\Bigl\{\sqrt{p\,|A|},\,
\frac{|A|^2}{p^{1/2}}\Bigr\}.
$$
In particular, if $|A|>p^{2/3},$ then
$$
|A(A+1)|\gg \sqrt{p\, |A|}.
$$
Let us show that this is optimal in general settings bound up to the implied
constant. Let $N<0.1p$ be a positive integer,
$M=[2\sqrt{Np}]$ and let $g$ be a generator of $F_p^*.$ Consider the
set
$$
X=\{g^{n}-1:\, n=1,2,\ldots, M\}.
$$
From the pigeon-hole principle, there is a number $L$ such that
$$
|X\cap\{g^{L+1},\ldots, g^{L+M}\}|\ge \frac{M^2}{2p}\ge N.
$$
Take
$$
A=X\cap\{g^{L+1},\ldots, g^{L+M}\}.
$$
Then we have $|A|\ge N$ and
$$
|A(A+1)|\le 2M\le 2\sqrt{pN}.
$$
Thus, it follows that for any positive integer $N<p$ there exists a
set $A\subset F_p$ with $|A|=N$ such that
$$
|A(A+1)|\ll\sqrt{p|A|}.
$$
This observation illustrates the precision  of our result for large
subsets of $F_p.$

When $|A|\cdot|B|\cdot|C|\approx p^2,$ Theorem~\ref{thm:optimal}
implies that
$$
|AB|\cdot|(A+1)C|\gg \sqrt{|A|^3\cdot |B|\cdot|C|}.
$$
This coincides with the bound that one can get when $A,B,C$ are
subsets of the set of real numbers $\mathbb{R}.$

\begin{theorem}
\label{thm:5/4} Let $A,B,C$ be finite subsets of \,
$\mathbb{R}\setminus\{0,-1\}.$ Then
$$
|AB|\cdot|(A+1)C|\gg \sqrt{|A|^3\cdot |B|\cdot |C|}.
$$
\end{theorem}
In particular, taking $B=A+1,\, C=A,$ we obtain the bound
$$
|A(A+1)| \gg |A|^{5/4}.
$$
We mention Elekes' sum-product estimate~\cite{El} in the case of
real numbers:
$$
|A+A|+|AA|\gg |A|^{5/4}.
$$
More generally Elekes' work implies that if $A,B,C$ are finite
subsets of the set
$\mathbb{R}\setminus\{0\},$ then
$$
|AB|\cdot |A+C|\gg \sqrt{|A|^3\cdot |B|\cdot|C|}.
$$
The best known bound up to date in the ``pure" sum-product problem
for real numbers is $|A+A|+|AA|\gg |A|^{4/3+o(1)},$ due to Solymosi~\cite{Sol}.

\section{Proof of Theorem~\ref{thm:106/105}}

For $E\subset A\times B$ we write
$$
A\En B=\{a-b: (a,b)\in E\}.
$$
A basic tool in the proof of Theorem~\ref{thm:106/105} is the
following explicit Balog-Szemer\'edi-Gowers type estimate given by
Bourgain and Garaev~\cite{BG}.

\begin{lemma}
\label{lem:BG1} Let $A\subset F_p, \, B\subset F_p,\, E\subset
A\times B$ be such that $|E|\ge |A||B|/K.$ There exists a subset
$A'\subset A$ such that $|A'|\ge 0.1 |A|/K$ and
$$
|A\En B|^4\ge \frac{|A'-A'|\cdot|A|\cdot|B|^{2}}{10^4K^{5}}.
$$
\end{lemma}

Theorem~\ref{thm:106/105} will be derived from the combination of
Lemma~\ref{lem:BG1} with the following specific variation of the
sum-product estimate from~\cite{BG}.
\begin{lemma}
\label{lem:BG2} Let $A\subset F_p,\, |A|<p^{1/2}.$ Then,
$$
|A-A|^8\cdot|A(A+1)|^4\ge |A|^{13+o(1)}
$$
\end{lemma}

The proof of Lemma~\ref{lem:BG2} follows from straightforward
modification of the proof of Theorem 1.1 of~\cite{BG}, so we only
sketch it. It suffices to show that
$$
|A-A|^5\cdot|2A-2A|\cdot|A(A+1)|^4\ge |A|^{11+o(1)}.
$$
Indeed, having this estimate established, one can apply it to large
subsets of $A,$ iterate the argument of Katz and Shen~\cite{KS}
several times and finish the proof; for more details, see~\cite{BG}.

We can assume that $A\cap \{0, -1\}=\emptyset$ and $|A|\ge 10.$
There exists a fixed element $b_0\in A$ such that
$$
\sum_{a\in A}|(a+1)A\cap (b_0+1)A|\ge \frac{|A|^3}{|A(A+1)|}.
$$
Decomposing into level sets, we get a positive integer $N$ and a
subset $A_1\subset A$ such that
\begin{equation}
\label{eqn:aAcapbAge} N\le |(a+1)A\cap (b_0+1)A|< 2N \quad {\rm for
\quad any} \quad a\in A_1,
\end{equation}
\begin{equation}
\label{eqn:N|A1|} N|A_1|\ge \frac{|A|^3}{2|A(A+1)|\cdot\log|A|}.
\end{equation}
In particular,
\begin{equation}
\label{eqn:boundA1} N\ge
\frac{|A|^2}{2|A|\cdot|A(A+1)|\cdot\log|A|}.
\end{equation}
We can assume that $|A_1|>1.$ Due to the observation of Glibichuk and Konyagin~\cite{GK}, either
$$
\frac{A_1-A_1}{A_1-A_1}=F_p
$$
or we can choose elements $b'_1,b'_2,b'_3,b'_4\in A_1$ such that
$$
\frac{b'_1-b'_2}{b'_3-b'_4}-1\not\in \frac{A_1-A_1}{A_1-A_1}.
$$
Using the step of Katz and Shen~\cite{KS}, we deduce that in either case there exist elements $b_1,b_2,b_3,b_4\in A_1$
such that
\begin{equation}
\label{eqn:length4} \Bigl|(b_1-b_2)A+(b_3-b_4)A\Bigr|\gg
\frac{|A_1|^3}{|A-A|}.
\end{equation}

To each element $x\in (b_1-b_2)A+(b_3-b_4)A$ we attach one fixed
representation
\begin{equation}
\label{eqn:attach x} x=(b_1-b_2)a(x)+(b_3-b_4)a'(x),\quad a(x),
a'(x)\in A.
\end{equation}
Denote
$$
S=(b_1-b_2)A+(b_3-b_4)A,\quad S_i=(b_i+1)A\cap (b_0+1)A; \quad
i=1,2,3,4.
$$
As in~\cite{BG}, we consider the mapping
$$
f: S\times S_1\times S_2 \times S_3\times S_4 \to (2A-2A)\times
(A-A)\times(A-A)\times(A-A)\times(A-A)
$$
defined as follows. Given
$$
x\in S, \quad x_i\in S_i; \quad i=1,2,3,4,
$$
we represent $x$ in the form~\eqref{eqn:attach x}, represent $x_i$
in the form
$$
x_i=(b_i+1)a_i(x_i)=(b_0+1)a_i'(x_i),\quad a_i(x_i)\in A,\quad
a_i'(x_i)\in A,\quad (i=1,2,3,4),
$$
and define
$$
f(x,x_1,x_2,x_3,x_4)=(u,u_1,u_2,u_3, u_4),
$$
where
$$
u=a_1'(x_1)-a_2'(x_2)+a_3'(x_3)-a_4'(x_4),
$$
$$
u_1=a(x)-a_1(x_1), \quad u_2=a(x)-a_2(x_2),
$$
$$
u_3=a'(x)-a_3(x_3),\quad u_4=a'(x)-a_4(x_4).
$$
From the construction we have
$$
x=(b_1+1)u_1-(b_2+1)u_2+(b_3+1)u_3-(b_4+1)u_4+(b_0+1)u.
$$
Therefore, the vector $(u,u_1,u_2,u_3,u_4)$ determines $x$ and thus
determines $a(x), a'(x)$ and consequently determines $a_1(x_1),
a_2(x_2), a_3(x_3), a_4(x_4)$ which determines $x_1,x_2,x_3,x_4.$
Hence, since $|(b_i+1)A\cap (b_0+1)A|\ge N,$ we get that
$$
|(b_1-b_2)A+(b_3-b_4)A|N^4\le |A-A|^4\cdot |2A-2A|.
$$
Taking into account~\eqref{eqn:length4}, we get
$$
|A-A|^4\cdot |2A-2A|\gg \frac{|A_1|^3N^4}{|A-A|}.
$$
Using~\eqref{eqn:aAcapbAge}--\eqref{eqn:boundA1}, we conclude the
proof of Lemma~\ref{lem:BG2}.

We proceed to prove Theorem~\ref{thm:106/105}. Denote
$$
E=\{(x, x+xy):\,\, x\in A, \, y\in A\}\subset A\times A(A+1),
$$
Then,
$$
|E|=|A|^2=\frac{|A|\cdot|A(A+1)|}{K},\quad K=\frac{|A(A+1)|}{|A|}.
$$
Let $B=A(A+1).$ Observe that
$$
-AA=A\En B.
$$

According to Lemma~\ref{lem:BG1} there exists $A'\subset A$ with
\begin{equation}
\label{eqn:A'end} |A'|\gg \frac{|A|}{K}=\frac{|A|^2}{|A(A+1)|}
\end{equation}
such that
$$
|AA|^4|A(A+1)|^3\gg |A'-A'||A|^6.
$$
Raising to eights power and multiplying by $|A(A+1)|^4\ge
|A'(A'+1)|^4,$ we get
$$
|AA|^{32}\cdot|A(A+1)|^{28}\gg |A'-A'|^8|A'(A'+1)|^4|A|^{48}.
$$
Combining this with Lemma~\ref{lem:BG2} (applied to $A'$), we obtain
$$
|AA|^{32}\cdot|A(A+1)|^{28}\gg |A'|^{13}|A|^{48+o(1)}.
$$
Taking into account the inequality~\eqref{eqn:A'end}, we get
$$
|AA|^{32}\cdot|A(A+1)|^{41}\ge |A|^{74+o(1)}.
$$
From Ruzsa's triangle inequalities in multiplicative form, we have
$$
|AA|\le\frac{|A(A+1)|\cdot|(A+1)A|}{|A+1|}=\frac{|A(A+1)|^2}{|A|}.
$$
Putting last two inequalities together, we conclude that
$$
|A(A+1)|^{105}\ge |A|^{106+o(1)}.
$$

\section{Proof of Theorem~\ref{thm:optimal}}

Let $J$ be the number of solutions of the equation
$$
x^{-1}y(z^{-1}t-1)=1, \quad (x,y,z,t)\in
AB\times B\times C\times (A+1)C.
$$
Observe that for any given triple $(a,b,c)\in A\times B\times C$ the
quadruple
$(x,y,z,t)=(ab, \, b, \, c, \, (a+1)c)$
is a solution of this equation. Thus,
\begin{equation}
\label{eqn:Jlower} J\ge |A|\cdot|B|\cdot|C|.
\end{equation}
On the other hand for any nonprincipal character $\chi$ modulo $p$
we have
$$
\Bigl|\sum_{z\in C}\,\,\sum_{t\in (A+1)C}\chi(z^{-1}t-1)\Bigr|\le
\sqrt{p\,|C|\cdot |(A+1)C|},
$$
see, for example, the solution to exercise 8 of~\cite[Chapter
V]{Vin}. Therefore, the method of solving multiplicative ternary
congruences implies that
$$
J=\frac{1}{p-1}\sum_{\chi}\sum_{x,y,z,t}\chi\Bigl(x^{-1}y(z^{-1}t-1)\Bigr)=
$$
$$
=\frac{1}{p-1}\sum_{x,y,z,t}\chi_0\Bigl(x^{-1}y(z^{-1}t-1)\Bigr)+
\frac{1}{p-1}\sum_{\chi\not=\chi_0}\sum_{x,y,z,t}\chi(x^{-1})\chi(y)\chi(z^{-1}t-1)
$$
$$
\le \frac{|AB|\cdot|B|\cdot |C|\cdot |(A+1)C|}{p-1}+\sqrt{p\,
 |C|\cdot|(A+1)C|\cdot|AB|\cdot |B|}.
$$
Comparing this with~\eqref{eqn:Jlower}, we conclude the proof.\\

{\bf Remark}. In Karatsuba's survey paper ~\cite{Kar} the interested reader will find many applications of character sums to multiplicative congruences.

\section{Proof of Theorem~\ref{thm:5/4}}

Since $A\cap \{0, -1\}=\emptyset,$ we can assume that $|A|$ is
large. We will use the Szemer\'edi-Trotter incidence theorem, which
claims that if $\mathcal{P}$ is a finite set of points $(x,y)\in
\mathbb{R}^2$ and $\mathcal{L}$ is a finite set of lines
$\ell\subset \mathbb{R}^2,$ then
$$
\#\Bigl\{\Bigl((x,y),\ell\Bigr)\in \mathcal{P}\times \mathcal{L}:\,
(x,y)\in \ell\Bigr\}\ll
|\mathcal{P}|+|\mathcal{L}|+(|\mathcal{P}||\mathcal{L}|)^{2/3}.
$$
We mention that this theorem was applied by Elekes in the above mentioned
work~\cite{El} to the sum-product problem for subsets of $\mathbb{R}.$ In application to our problem,
we let
$$
\mathcal{P}=\{(x,y):\, x\in AB,\, y\in (A+1)C\}
$$
and let $\mathcal{L}$ to be the family of lines $\{\ell=\ell(z,t):
z\in C,\, t\in B\}$ given by the equation
$$
y-\frac{z}{t}\,x-z=0.
$$
In particular,
$$
|\mathcal{P}|=|AB|\cdot|(A+1)C|,\quad |\mathcal{L}|=|B||C|.
$$ Each line $\ell(z,t)\in \mathcal{L}$ contains $|A|$ distinct points
$(x,y)\in \mathcal{P}$ of the form
$$
(x,y)=(at,\,(a+1)z);\quad a\in A.
$$
Thus,
$$
\#\Bigl\{\Bigl((x,y),\ell\Bigr)\in \mathcal{P}\times \mathcal{L}:\,
(x,y)\in \ell\Bigr\}\ge |A||\mathcal{L}|=|A|\cdot |B|\cdot |C|.
$$
Therefore, the Szemer\'edi-Trotter incidence theorem implies that
$$
|A|\cdot |B|\cdot |C|  \ll
|AB|\cdot|(A+1)C|+|B||C|+\Bigl(|AB|\cdot|(A+1)C|\cdot
|B|\cdot|C|\Bigr)^{2/3}.
$$
Since $|A|$ is large and  $|AB|\cdot|(A+1)C|\ge |A|^2,$ the result
follows.

\end{document}